\documentclass{article}

\usepackage[
  hyperindex  = {true},
  colorlinks  = {true},
  linkcolor   = {blue},
  citecolor   = {blue}
]{hyperref}
\usepackage{multirow}
\usepackage{color, colortbl}
\usepackage{caption}
\usepackage{subcaption}
\usepackage{latexsym,graphicx, amssymb, amsfonts, setspace, geometry, amsmath,amsthm}
\usepackage{cleveref}
\usepackage{pgfplots}
\pgfplotsset{compat=1.8}
\usepgfplotslibrary{statistics}
\usepackage{grffile}
\usepackage{pgfplotstable}
\usepackage{mathtools}
\usepackage[utf8]{inputenc}
\usepackage[english]{babel}
\usepackage{float}
\usepackage{esvect}
\usepackage{algorithm}
\usepackage{algpseudocode}
\usepackage{enumitem}
\usepackage{hhline}
\usepackage{varwidth}
\usepackage{tikz}
\usetikzlibrary{backgrounds,patterns,calc}
\usepackage{xparse}
\usepackage{tabularx}
\usepackage{placeins}

\pgfplotsset{compat=newest}

\makeatletter
\newcommand{\doublewidetilde}[1]{{%
  \mathpalette\double@widetilde{#1}%
}}
\newcommand{\double@widetilde}[2]{%
  \sbox\z@{$\m@th#1\widetilde{#2}$}%
  \ht\z@=.9\ht\z@
  \widetilde{\box\z@}%
}
\makeatother

\newtheorem{remark}{Remark}[section]
\newtheorem{theorem}{Theorem}[section]

\newcommand{\rb}{\boldsymbol{r}}
\newcommand{\xb}{\boldsymbol{x}}
\newcommand{\ub}{\boldsymbol{u}}

\newcommand{\tran}{^{\top\kern-\scriptspace}}
\newcommand{\bnu}{\boldsymbol{\nu}}

\newcommand{\be}{\begin{equation}}
\newcommand{\ee}{\end{equation}}
\newcommand{\ba}{\begin{aligned}}
\newcommand{\ea}{\end{aligned}}
\newcommand{\bea}{\begin{eqnarray}}
\newcommand{\eea}{\end{eqnarray}}

\newcommand\nc{\newcommand}

\newcommand\niter{n_{\textrm{iter}}}
\newcommand\nitmax{n_{\textrm{max-iter}}}

\newcommand{\btheta}{\boldsymbol{\theta}}

\newcommand{\bE}{\mathbf{E}}
\newcommand{\bH}{\mathbf{H}}

\newcommand{\uinc}{u^\emph{inc}}

\nc\ex{E_x}
\nc\ey{E_y}
\nc\ez{E_z}

\nc\hx{H_x}
\nc\hy{H_y}
\nc\hz{H_z}

\nc\px[1]{\frac{\partial #1}{\partial x}}
\nc\py[1]{\frac{\partial #1}{\partial y}}

\nc\br{{\boldsymbol{r}}}
\nc\curl{\nabla\times}
\nc\dive{\nabla\cdot}
\nc\pt{\frac{\partial}{\partial t}}
\nc\pet{\frac{\partial \bE}{\partial t}}
\nc\pht{\frac{\partial \bH}{\partial t}}

\newcommand{\Lim}[1]{\raisebox{0.5ex}{\scalebox{0.8}{$\displaystyle \lim_{#1}\;$}}}

\newcommand{\bx}{\boldsymbol{x}}

\definecolor{Gray}{gray}{0.9}
\newcolumntype{g}{>{\columncolor{Gray}}c}

\pgfmathdeclarefunction{gauss}{3}{%
  \pgfmathparse{1/(#3*sqrt(2*pi))*exp(-((#1-#2)^2)/(2*#3^2))}%
}

\title{On the robustness of inverse scattering for penetrable, homogeneous objects with complicated
boundary}

\author{Carlos Borges\thanks{Department of Mathematics, University of Central Florida, Orlando, FL, USA. \textit{Email: carlos.borges@ucf.edu}}
\and
Manas Rachh\thanks{Flatiron Institute, New York, NY, USA. \textit{Email: mrachh@flatironinstitute.org}}
\and
Leslie Greengard\thanks{Courant Institute of Mathematical Sciences, New York University, New York, NY, USA and Flatiron Institute, New York, NY, USA. \textit{Email: greengard@cims.nyu.edu}}
}

\begin{document}

\maketitle

\begin{abstract}
The acoustic inverse obstacle scattering problem consists of determining 
the shape of a domain from measurements of the scattered far field due to 
some set of incident fields (probes).
For a penetrable object with known sound speed, this can be accomplished by 
treating the boundary alone as an unknown curve.  Alternatively, one can 
treat the entire object as unknown and use a more general {\em volumetric}
representation, without making use of the known sound speed.
Both lead to strongly nonlinear and nonconvex optimization problems for which
recursive linearization provides a useful framework for numerical analysis.
After extending our shape optimization approach developed earlier for 
impenetrable bodies, we carry out a systematic study of both methods and 
compare their performance on a variety of examples. Our findings indicate 
that the volumetric approach is more robust, even though the 
number of degrees of freedom is significantly larger. 
We conclude with a discussion of this 
phenomenon and potential directions for further research.
\end{abstract}


\section{Introduction}
Using waves as probes for non-destructive or non-invasive testing is of interest in a diverse 
set of applications, from medical imaging to materials characterization, remote sensing, sonar
and radar \cite{kuchment2014radon,collins1995nondestructive,engl2012inverse,Ustinov2014,cheney2009fundamentals,
beilina2014reconstruction,klibanov2019convexification,thanh2014reconstruction}. 
With acoustic waves, one can imagine recovering the shape, density, and/or internal sound speed 
from measurements of the scattered field induced by a collection of known incident waves.
In this paper, we focus on recovering the shape
$\Gamma$ of a penetrable obstacle $\Omega_i$ with known density $\rho_i$ and sound speed $c_i$ 
immersed in a medium $\Omega_e=\mathbb{R}^2\setminus\Omega_i$ with density $\rho_e$ and sound speed $c_e$, assuming the data consists of far field measurements when a plane wave impinges on the
unknown obstacle from multiple directions at multiple frequencies. 
We will refer to this as the {\em inverse penetrable obstacle problem} or 
simply the {\em inverse obstacle problem}.
The incident wave at frequency $\omega$ will be denoted by $u^\emph{inc}$, the field
in the interior will be denoted by $u_i$ and the scattered field will be denoted by $u_e$.
The total field $u$ is equal to $u_i$ in the interior and to the sum $u^\emph{inc} + u_e$ in the 
exterior.
The corresponding (time harmonic) forward problem at frequency $\omega$
consists of solving
\begin{equation} \label{eq:transmission_problem}
\begin{cases}
\Delta u_i+k_i^2u_i = 0,  \quad \text{in} ~\Omega_i, \\
\Delta u_e+k^2u_e = 0, \quad \text{in} ~\Omega_e, \\ 
\left[u\right] = 0 ~\text{and} ~\left[\beta \frac{\partial u}{\partial \nu} \right] = 0, \quad \text{on}~ \Gamma, \\
\Lim{\|{\bf x}\|\rightarrow \infty} \|{\bf x}\|^{1/2}\left(\frac{\partial u_e}{\partial r} - ik u_e\right) = 0.
\end{cases} 
\end{equation}
assuming the density and sound speed are known.
Here, $k=\omega/c_e$ is the wavenumber in the exterior medium and $k_i=\omega/c_i$ 
is the wavenumber in the interior of the obstacle,
$\left[u\right]=u\vert_+-u\vert_-$, with $u\vert_+$ and $u\vert_-$ denoting the limit of the 
total field as one approaches $\Gamma$ from the exterior and interior, 
respectively. $\beta$ is a piecewise constant function with $\beta=1$ in $\Omega_e$ 
and $\beta=\rho_e/\rho_i$ in $\Omega_i$. 
The direction of the incident plane wave will be denoted by the unit vector $\btheta$, so that
$u^\emph{inc}(\xb)=e^{i k\xb\cdot\btheta}$.
The solution to \eqref{eq:transmission_problem} can be obtained using a boundary
integral equation method \cite{hohage1998newton, hettlich1995frechet, greengard2012stable, colton2013integral, rokhlin1983solution} with two unknown source densities
supported on $\Gamma$ alone.
The obstacle problem and the integral equation approach are described in section 
\ref{s:trans_problem}.

Alternatively, one can treat the domain as having an unknown shape defined by the (compact) support
of some perturbation, $q(x)$, of the sound speed in the ambient space. We will refer to this as
the {\em inverse medium problem} \cite{colton2019inverse}.
The forward problem here consists of computing the scattered field $u^\emph{scat}$ induced by the
same incident wave $u^\emph{inc}$ as above when it impinges on the obstacle defined by 
a {\em known} function $q(x)$:
\begin{equation} \label{eq:volume_problem}
\begin{cases}
\Delta u^\emph{scat}+k^2(1+q)u^\emph{scat} = -k^2 q u^\emph{inc},  \quad \text{in} ~\mathbb{R}^2, \\
\Lim{\|{\bf x}\|\rightarrow \infty} \|{\bf x}\|^{1/2}\left(\frac{\partial u^\emph{scat}}{\partial r} - iku^\emph{scat}\right) = 0.
\end{cases} 
\end{equation}
While problem \eqref{eq:volume_problem}, like 
problem \eqref{eq:transmission_problem}, can be solved using integral equation techniques,
the scattered field must now be represented as a volume potential, leading to the 
Lippmann-Schwinger integral equation \cite{colton2019inverse}.
This approach is discussed in section \ref{s:vol_problem}.

\begin{remark}
If $\rho_i = \rho_e$ in problem \eqref{eq:transmission_problem},
the support of the function $q$ is 
defined to be $\Omega_i$, and $k_i^2=k^2(1+q)$ inside $\Omega_i$, then the two forward
problems are identical.
\end{remark}

Suppose now that a collection of $N_r$ receivers are equispaced on a disk of radius $R \gg 1$: 
$\{ \br_{m} = (R \cos \frac{2\pi m}{N_r}, R \sin \frac{2 \pi m}{N_r}), \  m=1,2\ldots N_r \}$.
We then define the forward transmission operator by
\begin{equation} \label{eq:fwd_transm_op}
\mathcal{F}^{(T)}_{k,\btheta}(\Gamma)=\ub_{k,\btheta}^\emph{meas} \in  \mathbb{C}^{N_r} ,
\end{equation}
where the $m$th component of $\ub_{k,\btheta}^\emph{meas}$ is the scattered field $u_e(\br_m)$,
that solves~\eqref{eq:transmission_problem} in response to
an incoming plane wave $u^\emph{inc}(\bx)=e^{ik\xb\cdot\btheta}$.
The inverse obstacle problem can now be stated more precisely in terms of optimization,
namely in the form
\begin{equation}
\tilde{\Gamma}=\arg \min_{\Gamma} \sum_{k_l,\btheta_j} \| \ub_{k_l,\btheta_j}^{\emph{meas}} -\mathcal{F}^{(T)}_{k_l,\btheta_j}(\Gamma)\|, \label{eq:transmission_inv_problem}
\end{equation}
assuming we have $N_k$ probing frequencies $\{k_1,\dots,k_{N_k} \}$
and $N_d$ incident directions $\{\btheta_1,\dots,\btheta_{N_k} \}$.

Similarly, the forward volumetric scattering operator is defined by
\begin{equation} \label{eq:fwd_vol_op}
\mathcal{F}^{(V)}_{k,\btheta}(q)=\ub_{k,\btheta}^\emph{meas} \in  \mathbb{C}^{N_r} ,
\end{equation}
where the $m$th component of $\ub_{k,\btheta}^\emph{meas}$ is the scattered field 
$u^\emph{scat}(\br_m)$ that solves ~\eqref{eq:volume_problem},
in response to an incoming plane wave
$u^\emph{inc}(\bx)=e^{ik\xb \cdot \btheta}$ impinges on the inhomogeneity. 
The inverse medium problem can be stated in the form
\begin{equation}
\tilde{q}=\arg \min_{q} \sum_{k_l,\btheta_j} \| \ub_{k_l,\btheta_j}^{\emph{meas}} -\mathcal{F}^{(V)}_{k_l,\btheta_j}(q)\|.\label{eq:volume_inv_problem}
\end{equation}

The inverse problems \eqref{eq:transmission_inv_problem} and \eqref{eq:volume_inv_problem} are 
fully nonlinear, ill-posed without additional constraints, and computationally challenging. 
A variety of nonlinear iterations have been applied to such problem, including the
Gauss-Newton method \cite{colton2019inverse}, other Newton-like variants, Landweber iteration, and
the nonlinear conjugate gradient method \cite{hohage,kress2007}.
To deal with the ill-posedness, regularization methods such as 
Tikhonov regularization or the truncated SVD can be used. Alternatively, one can use a 
parametric approximation of the object with fewer degrees of freedom in the forward model. One such
approximation is to assume that the object is bandlimited, with the bandlimit determined by the 
frequency and the number of independent measurements.
Finally, one can either choose to solve for the unknown object using all frequencies 
simultaneously or solve a sequence of single frequency inverse problems. 
Since using all frequencies together is computationally expensive, Chen 
\cite{chen1995recursive, chen1997inverse} suggested the recursive linearization algorithm (RLA).
In this approach, 
one first solves for a low-resolution approximation of the unknown using only the
lowest frequency data. This reconstruction then serves as an initial guess for Gauss-Newton iteration at the next available frequency, until the highest frequency data has been reached. At each step
of this iteration, the complexity of the unknown is gradually increased, typically by increasing its
bandlimit. The RLA has been successfully applied 
to both inverse obstacle scattering and the inverse medium problem
\cite{bao2007inverse, bao2005inverse, bao2010error, bao2011imaging, bao2010multi, bao2011numerical,
 bao2015recursive, borges2017high, chaillat2012faims, borges2015inverse, sini2012inverse, 
borges2022multifrequency}.  We refer the reader to \cite{bao2015inverse} for
a thorough review of inverse scattering problems based on multiple frequency data. 

The present paper is aimed at a question that appears not to have been considered previously:
namely, to compare inverse obstacle scattering and the inverse medium problem as numerical
approaches when the forward problems are identical, as outlined above.
In section \ref{s:trans_problem}, after reviewing integral equation methods for forward scattering
from a penetrable obstacle, we extend the recursive linearization method of 
\cite{sini2012inverse, borges2015inverse,borges2022multifrequency} for sound-soft
obstacles to the current setting.
Next, we briefly review the volumetric inverse scattering method (also based on recursive
linearization) presented in \cite{borges2017high}.
Both methods have been shown to be capable of obtaining high resolution reconstructions of very 
complicated, but suitably bandlimited, unknowns. 
Focusing on computational complexity, it is 
straightforward to see that the inverse obstacle approach should be much faster. Boundary
integral methods are used for the forward scattering problem and only an unknown curve is being
sought. The volumetric inverse scattering approach requires volume integral equations to be 
solved and seeks an unknown function defined on a two-dimensional region in the plane.

Note, however, that the formulations have distinct features, ignoring questions of computational
efficiency. First, the obstacle scattering approach is seeking a discontinuity in sound speed
defined on a smooth curve, while the volumetric scattering approach is seeking a smooth
(bandlimited) function with the same scattered far field.
It is easy to imagine that this affects the robustness of the solver, 
the dependence on the number of measurements, etc. 
We will explore these questions numerically in 
section \ref{s:num_res}. 
An interesting discovery is that there are clear cases where the 
volumetric approach is able to obtain high quality reconstruction while the obstacle approach fails.
We conclude with a discussion of these results and opportunities for further research in
section \ref{s:conclusions}.

\begin{remark} 
Since the data sets are somewhat complicated, and depend on several 
parameters, we summarize some of the important notation
in Table \ref{table:symbol}. 
When the context is clear, we will omit some indices and write, for example,
$\mathcal{F}^{(V)}_k(q) \coloneqq 
\left[\mathcal{F}^{(V)}_{k,\btheta_1}(q); \dots;\mathcal{F}^{(V)}_{k,\btheta_{N_d}}(q)\right]$ 
or ${\bf u}^{\emph{meas}}$ to refer to the vector whose $m$th 
component is ${\bf u}^{\emph{meas}}_{k,\btheta}(\br_m)$ for incident direction $\btheta$ at
wavenumber $k$.

\begin{table}
\caption{List of principal symbols used in this article.}\label{table:symbol}
{\small
\begin{center}
\begin{tabular}{ll}
\hline
Symbol & Description \\
\hline\hline
$\Omega_i$          & interior of the penetrable obstacle\\
$\Omega_e$          & exterior of the penetrable obstacle\\
$\Gamma$            & boundary of $\Omega_i$ \\
$q$                 & perturbation of the sound speed in $\Omega_i$ \\
$\omega$            & frequency of the incident plane wave \\
$k$               & wavenumber of the incident plane wave in $\Omega_e$ \\
$k_i$               & wavenumber of the incident plane wave in $\Omega_i$\\
$\rb_j$         & location of $j$th  receiver in $\Omega_e$\\
$\btheta$       & incident direction of plane wave $u^{\emph{inc}}_{k,\btheta}$ ($\|\btheta\|=1$) \\
$u^{\emph{inc}}_{k,\btheta}$  &incident plane wave with wavenumber $k$ and direction $\btheta$\\
$u^{\emph{scat}}_{k,\btheta}$   &scattered field generated by $u_{k,\btheta}^{\emph{inc}}$ \\
${\bf u}^{\emph{meas}}_{k,\btheta}$  &vector in $\mathbb{C}^{N_r}$, with $j$th
component $u^{\emph{scat}}_{k,\btheta}(\br_j)$ \\
$N_k$                   &number of probing frequencies\\
$N_r$			&number of receivers \\
$N_d$                   &number of incident waves \\
$\mathcal{F}^{(V)}_{k,\btheta}$  &forward operator for volume scattering mapping $q$ in $\Omega_i$ to ${\bf u}^{\emph{meas}}_{k,\btheta}$  (for given $u_{k,\btheta}^{\emph{inc}}$) \\
$\mathcal{F}^{(T)}_{k,\btheta}$  &forward operator for obstacle scattering mapping $\partial\Omega$ to ${\bf u}^{\emph{meas}}_{k,\btheta}$ (for given $u_{k,\btheta}^{\emph{inc}}$) \\
$\mathcal{J}^{(V)}_{k,\btheta}$ &Frech\'{e}t derivative of $\mathcal{F}^{(V)}_{k,\btheta}$ with respect to the function $q$ \\
$\mathcal{J}^{(T)}_{k,\btheta}$ &Frech\'{e}t derivative of $\mathcal{F}^{(T)}_{k,\btheta}$ with respect to the boundary $\Gamma$ \\
$\mathcal{S}$         & single layer potential  \\
$\mathcal{D}$         & double layer potential \\
$\mathcal{K}$         & normal derivative of the single layer potential  \\
$\mathcal{T}$         & normal derivative of the double layer potential  \\
$G$                   & free space Green's function for the two dimensional Helmholtz equation \\
$N_{\gamma,k}$	      & number of Fourier modes used to update the curve $\Gamma$ at wavenumber $k$\\
$N_m$                 & number of modes in sine series used to update the function $q$ \\

\hline 
\end{tabular}
\end{center}
}
\end{table}

\end{remark}

\section{Inverse penetrable obstacle scattering}\label{s:trans_problem}

In this section, we first describe an integral equation approach to solving the forward problem
(under the assumption that $\rho_i = \rho_e$).
We then review the recursive linearization method for the inverse problem (shape recovery) using
multifrequency data. 
The forward transmission scattering operator, defined in the introduction, generates the data
vector ${\bf u}_{k,\btheta}^\emph{meas} = (u_e({\rb}_1),u_e({\rb}_2),\dots,u_e({\rb}_{N_r}))$
for an incident plane wave $u^{\emph{inc}}_{k,\btheta} = e^{ik \bx \cdot \btheta}$. 
This requires solving problem \eqref{eq:transmission_problem}, for which we use potential theory.
We represent the exterior and interior fields as
$u_e=\left(\mathcal{D}_k \sigma - \mathcal{S}_k \mu\right)$ and 
$u_i=  \left( \mathcal{D}_{k_i} \sigma - \mathcal{S}_{k_i} \mu \right)$, where 
\begin{equation}
\mathcal{D}_{k} \sigma = \int_{\Gamma} \frac{\partial G(k \|{\bf x}-{\bf y}\|)}{\partial \nu(y)} \sigma({\bf y})ds(y),\quad\text{and}\quad
\mathcal{S}_{k} \mu = \int_{\Gamma} G(k \|{\bf x}-{\bf y}\|) \mu({\bf y})ds(y),
\end{equation}
are the single and double-layer potentials
for the domain with wavenumber $k$. Here, $G(k r)=H_0^{(1)}(k r)$ is 
the Green's function for the two dimensional Helmholtz equation in free-space, where
$H_0^{(1)}$ denotes the first kind Hankel function of order zero satisfying the Sommerfeld 
radiation condition at infinity.
Enforcing the interface conditions on $\Gamma$ leads to the system of linear integral equations
\begin{equation} \label{eq:transmission_eq}
\begin{split}
\left[ \mathcal{I}+\left(
 \mathcal{D}_k-  \mathcal{D}_{k_i} \right)\right]\sigma+
\left( \mathcal{S}_{k_i}-  \mathcal{S}_k \right) \mu = - \uinc\\
\left(\mathcal{D}_k^\prime-\mathcal{D}_{k_i}^\prime\right)\sigma+
\left(\mathcal{I}+\mathcal{S}_{k_i}^\prime-\mathcal{S}_k^\prime\right)\mu = 
- \frac{\partial \uinc}{\partial \nu},
\end{split}
\end{equation}
where $\mathcal{I}$ is the identity operator, 
$\mathcal{D}_{k}^\prime$ and $\mathcal{S}_{k}^\prime$ 
are the normal derivatives of the operators $\mathcal{D}_{k}$ and 
$\mathcal{S}_{k}$ and all operators 
in \eqref{eq:transmission_eq} are interpreted in a principal value sense.

The system above is a Fredholm integral equation of the second kind for the 
unknown densities $(\sigma, \mu)$. It is well-known to have a unique solution
\cite{colton2013integral}.
For numerical purposes, we discretize the boundary at equispaced points along $\Gamma$, using
Alpert's 16$th$ order Gauss-trapezoidal rule \cite{alpert1999hybrid} for quadrature.
To ensure more than 10 digits of accuracy, it suffices to use 10 points per wavelength or more, with the wavelength taken 
to correspond to the larger of the two wavenumbers 
$k_i$ and $k$. In all our examples, we used 70 points per wavelength, unless stated otherwise. 
For simplicity, we solve the resulting system of equations using standard 
LU factorization, so that the complexity of solution is of the order $\mathcal{O}(N^3)$, 
where $N$ is the total number of points along the boundary. If $N$ were larger than in the 
examples considered here, one could replace Gaussian elimination with a fast direct solver that has
$O(N \log N)$ complexity
\cite{bebendorf2005hierarchical, borm2003hierarchical, borm2003introduction, chandrasekaran2006fast, chandrasekaran2006fast1, gillman2012direct, greengard2009fast, ho2012fast, martinsson2005fast}. 
That approach is discussed for large-scale inverse scattering problems with sound-soft obstacles in \cite{borges2015inverse}. 

\subsection{Inverse obstacle scattering}

To recover the shape of the boundary of $\Gamma$, we now consider the optimization 
problem \eqref{eq:transmission_inv_problem}. To avoid the expense and complexity
of the
full multifrequency problem, as noted earlier, we solve a sequence of single frequency
problems using recursive linearization \cite{chen1995recursive, chen1997inverse}.
Thus, at each stage, corresponding to exterior wavenumber $k_m$, we must solve
the nonlinear, nonconvex, ill-posed problem
\begin{equation}
\tilde{\Gamma}=\arg \min_{\Gamma} \sum_{\theta_j=1}^{N_d} \| \ub_{k_m,\theta_j}^{\emph{meas}} -\mathcal{F}^{(T)}_{k_m,\theta_j}(\Gamma)\|. \label{eq:transmission_single_inv_problem}
\end{equation}
Assuming that we have solved the problem correctly at frequency $k_{m-1}$, 
the obtained approximation $\Gamma^{(m-1)}$
to $\Gamma$ should 
provide a good initial guess for the unknown curve $\Gamma^{(m)}$.
That is, we are using frequency as a homotopy parameter, typically using 
steepest descent or Gauss-Newton iteration
to solve for the update $\delta \Gamma^{(m)} = \Gamma^{(m)} - \Gamma^{(m-1)}$. 
For sufficiently small steps in frequency, this overcomes
the intrinsic nonconvexity of \eqref{eq:transmission_single_inv_problem}
(although without a guarantee of global convergence).
To address the ill-posedness of the problem, as we march in frequency, we gradually increase
the complexity of the curve, parameterized in arclength.

More concretely, at the first stage, using the lowest frequency $k_0$,
and with no prior information about the curve, we let our initial guess be the
unit disk and seek a perturbation $h_1(t)$ in the normal direction to better fit the data. 
That is, we write 
\[ \Gamma^{(0)}(t) = e^{it}, \quad
h_1(t) =  \alpha_0 + \alpha_1 \cos(t) + \beta_1 \sin(t),\quad
\delta \Gamma^{(1)}(t) \equiv h_1(t) \bnu^{(0)}(t).
 \]
for $t \in [0,2\pi]$,
With a slight abuse of notation,
we make use of the equivalence of $\mathbb{R}^2$ and the complex plane, 
and view $\Gamma^{(0)}(t)$ as a complex-valued function whose real and imaginary parts
are the $x$ and $y$ components of the curve. 
The vector $\bnu^{(0)}(t)$ denotes the outward normal to $\Gamma^{(0)}(t)$, also 
treated as a vector in the complex plane.
Ignoring for the moment how we actually find 
$h_{1}(t)$, we define  
$\Gamma^{(1)}(t) = \Gamma^{(0)}(t) + \delta \Gamma^{(1)}(t)$ and determine its arclength
$L_1$. We then reparametrize the curve as 
\[ \Gamma^{(1)}(t) = 
\sum_{n=-N_{\gamma,k_1}}^{N_{\gamma,k_1}} \gamma^{(1)}_n e^{2 \pi int/L_1} .
 \]
Here, $N_{\gamma,k_1} = \max(1, \frac{p L_1 k_1}{2\pi})$, where $p$ denotes the desired number of
points per wavelength (which is typically in the range 10 to 70, as noted above).
The general step follows naturally. Since we may do more than one iteration at some wavenumber
$k_m$, let us denote by $j$ the iteration number,
by $\Gamma^{(j-1)}(t)$ the previous approximation to the curve with normal
$\bnu^{(j-1)}(t)$,
and by
\begin{equation}
 h_{j}(t) = \alpha_0^{(j)} + \sum_{n=1}^{\lceil c k_m \rceil} \alpha_n^{(j)} \cos(2\pi n t/L_{j-1})
+ \beta_n^{(j)} \sin(2\pi n t/L_{j-1})
\label{eq:representation_step_tranmission}
\end{equation}
the newly computed perturbation in the normal direction (by a method to be discussed shortly).
We typically choose $c=2$ and recommend keeping $c \in [1,3]$.
After obtaining the new arclength $L_j$, one sets
\[ 
\delta \Gamma^{(j)}(t) = h_j(t) \bnu^{(j-1)}(t),
 \]
 and reparametrizes the curve as a Fourier series in arclength as
\begin{equation}
 \Gamma^{(j)}(t) = \Gamma^{(j-1)}(t) + \delta \Gamma^{(j)}(t) \approx
\sum_{n=-N_{\gamma,k_m}}^{N_{\gamma,k_m}} \gamma^{(j)}_n e^{2 \pi int/L_j} ,
\label{eq:curveparam}
\end{equation}
where 
\begin{equation}
N_{\gamma,k_m} = \max(N_{\gamma,k_{m-1}}, \frac{p L_j k_m}{2\pi}).
\label{eq:nmodes}
\end{equation}

To further stabilize our nonlinear search, we augment the unconstrained formulation
\eqref{eq:transmission_single_inv_problem} with a trust region. For this, with 
$\Gamma$ a smooth curve of length $L$ with curvature
\[ \kappa(t) = \sum_{n=-\infty}^{\infty} \kappa_n e^{2 \pi i n t/L},
\]
we define 
\[ {\cal E}_\kappa \equiv \sum_{n=-\infty}^{\infty} |\kappa_n|^2 , \quad
{\cal E}_\kappa^{M} \equiv 
\sum_{n=-{M}}^{{M}} 
 |\kappa_n|^2.
\]
We will refer to ${\cal E}_\kappa$ as the elastic energy in the curve 
and
${\cal E}_\kappa^{M}$ as the elastic energy of its band-limited approximation.
We define
the set of allowed (closed, non-self-intersecting) curves at wavenumber $k$ (and parameter $c$) by 
\begin{equation}
{\cal A}_{k,\epsilon_f} =
\{ \Gamma\ |\ 
{\cal E}_\kappa^{(\lceil c k \rceil)} \geq (1-\epsilon_f) {\cal E}_\kappa \},
\label{eq:trust}
\end{equation}
i.e., ${\cal A}$ is the set of simple closed 
curves at wavenumber $k$ whose bandlimited approximation
captures the bulk $(1-\epsilon_f)$ of the elastic energy of the curve.
\begin{remark} \label{rmk:heisenberg}
Keeping the number of degrees of freedom in $h_{j}(t)$ and $\Gamma^{(j)}(t)$ proportional to 
$k_m$ as one increases the
incident frequency, and restricting the elastic energy of the curve to nearly bandlimited are forms of regularizations. 
This mitigates the ill-posedness
of \eqref{eq:transmission_single_inv_problem} in a physically sensible manner,
since the signature of high frequency features of the geometry (those that exceed the
frequency of the incident wave) decay exponentially in the far field
(the Heisenberg uncertainty principle for waves). Attempting to recover those features
is unstable.
\end{remark}

We return now to the optimization problem itself and write
the linearization of the system
\begin{equation}
\mathcal{F}^{(T)}_{k,\theta_j}(\Gamma^{(j-1)} + \delta \Gamma^{(j)}) = \ub_{k,\theta_j}^{\emph{meas}}
\label{eq:tsystem}
\end{equation}
as 
\begin{equation}
\mathcal{J}^{(T)}_{k,\theta_j} \, \delta \Gamma^{(j)} = \ub_{k,\theta_j}^{\emph{meas}}
- \mathcal{F}^{(T)}_{k,\theta_j}(\Gamma^{(j-1)}),
\label{eq:tsystem}
\end{equation}
where $\mathcal{J}^{(T)}_{k,\theta_j}$ is the Fr\'{e}chet derivative of the operator
$\mathcal{F}^{(T)}_{k,\theta_j}$ for the current guess $\Gamma^{(j-1)}$.
For a single angle of incidence, this is an underdetermined system. Using all incident angles,
we have the nonlinear least squares problem
\begin{equation} \label{eq:nonlin}
\mathcal{J}^{(T)}_k \delta \Gamma = \ub_{k}^{\emph{meas}} - \mathcal{F}_k^{(T)}(\Gamma^{(j-1)}),
\end{equation}
where 
\begin{align*}
 \mathcal{J}^{(T)}_k &=\left(\mathcal{J}^{(T)}_{k,\theta_1}; \cdots; \mathcal{J}^{(T)}_{k,\theta_{N_d}}\right),  \\
\ub_{k}^{\emph{meas}} &=\left(\ub_{k,\theta_1}^{\emph{meas}};\cdots;\ub_{k,\theta_{N_d}}^{\emph{meas}}\right), \\
\mathcal{F}_k^{(T)}(\Gamma^{(j-1)}) &=\left(\mathcal{F}_{k,\theta_1}^{(T)}(\Gamma^{(j-1)}); \cdots; \mathcal{F}_{k,\theta_{N_d}}^{(T)}(\Gamma^{(j-1)})\right).
\end{align*}

Following the discussion of \cite{powell1970new, powell1970hybrid},
the Gauss-Newton solution to \eqref{eq:nonlin} is given by
\begin{equation} \label{eq:gn_step}
\delta \Gamma_{GN} = 
\left(   (\mathcal{J}^{(T)}_k)^* 
\mathcal{J}^{(T)}_k \right)^{-1}  
({\mathcal{J}^{(T)}}_k)^*
\left(\ub_{k}^{\emph{meas}} - \mathcal{F}_k^{(T)}(\Gamma^{(j-1)}) \right)
\end{equation}
and the steepest descent direction by
\begin{equation} \label{eq:sd_step}
\delta \Gamma_{SD} = \left(\mathcal{J}^{(T)}_k\right)^* \left(\ub_{k}^{\emph{meas}} - \mathcal{F}_k^{(T)}(\Gamma^{(i)})\right).
\end{equation}

Rather than use the Gauss-Newton solution in an unconstrained fashion,
as in \cite{borges2015inverse,borges2022multifrequency},
we modify Powell's dogleg method \cite{powell1970new, powell1970hybrid}, which makes use of a trust 
region and both the Gauss-Newton and steepest descent steps.
In our proposed method, we first calculate $\delta \Gamma_{GN}$ and $\delta \Gamma_{SD}$. 
We then define $\Gamma^{(j)}_{GN}=\Gamma^{(j-1)}+\delta\Gamma_{GN}$ and 
$\Gamma^{(j)}_{SD}=\Gamma^{(j-1)}+\delta \Gamma_{SD}$ and check 
the elastic energies ${\cal E}_\kappa, {\cal E}_\kappa^{(\lceil c k \rceil)}$ 
of the two curves to determine if they lie in the trust region
${\cal A}_{k,\epsilon_f}$. We proceed as follows:

\begin{itemize}
\item If both curves lie in the trust region, we calculate the residual
$\sum_{\theta_j=1}^{N_d} \|\mathcal{F}_{k,\theta_j}^{(T)}(\Gamma)-\ub_{k,\theta_j}^{\emph{meas}}\|$ 
for both $\Gamma^{(j)}_{GN}$ and $\Gamma^{(j)}_{SD}$ and choose the step with the smaller residual.
\item If only one of the updated curves lies in the trust region, we accept that step. 
\item If neither curve lies in the trust region, we apply a Gaussian filter to the 
update $\delta \Gamma$. If one or both of the filtered curves
$\Gamma^{(j)}_{GN}$ and $\Gamma^{(j)}_{SD}$ lie in the trust region, we continue as above. 
Otherwise, we repeat the filtering up to 
$10$ times.
The Gaussian filter is defined by 
\begin{equation}
\gamma^{(j)}_{n} \to \gamma^{(j)}_{n} exp^{-\frac{n^2}{\sigma^2 N_{\gamma, k_{m}}^2}} \, ,
\end{equation}
where $\gamma^{(j)}_{n}$ are the coefficients defining the update $\delta \Gamma$, and $\sigma = 1/10^{(\ell-1)}$, 
where $\ell$ is the iteration number for the filtering step.
\end{itemize}

\begin{remark}
The filtering step here
is consistent with the discussion in Remark \ref{rmk:heisenberg}. If the elastic energy of the curve
is not captured by the first $N_{\gamma,k_m}$ modes, it must have a nontrivial evanescent far-field
signature and we are seeking the most band-limited curve that accurately reproduces the measurements.
Thus, filtering is consistent with our search regularization strategy.
See \cite{borges2020inverse, borges2022multifrequency} for further discussion of this point
and a more complete description of the algorithm.
\end{remark}

Finally, we conclude this section with a theorem that explains how one actually computes the 
Fr\'echet derivative $\mathcal{J}_{k,\theta_j}^{(T)}$ of the operator 
$\mathcal{F}_{k,\theta_j}^{(T)}$ for a given curve $\Gamma$. 
More precisely, we state how to compute the action of 
$\mathcal{J}_{k,\theta_j}^{(T)}$ on a normal perturbation $h(t) \bnu(t)$.
The proof follows very closely that presented in 
\cite{hettlich1995frechet, hohage1998newton}.

\begin{theorem}
Assume that $u_i$ and $u_e$ are the solutions to the transmission problem 
\eqref{eq:transmission_problem}, with incoming field $u^\emph{inc}(\xb)=e^{ik \xb\cdot\theta}$,
and let $\delta \Gamma = h(s) \bnu(s)$ denote a perturbation to the smooth 
curve $\Gamma$ in the normal
direction.
Then, the operator $\mathcal{F}^{(T)}_{k,\theta}$ is 
Fr\'echet differentiable at $\Gamma$, and the product of its Jacobian with a normal
perturbation $\delta \Gamma$, 
$\mathcal{J}^{(T)}_{k,\theta} \, \delta\Gamma$, is given at the 
receiver locations by the solution to the following boundary value problem: 
\begin{equation}\label{eq:der_transmission}
\begin{cases}
\Delta v_i+k_i^2v_i = 0,  \quad \text{in} ~\Omega_i, \\
\Delta v_e+k^2v_e = 0, \quad \text{in} ~\Omega_e, \\ 
v_{i-}-v_{e+}= h \, \left(\frac{\partial u_{+}}{\partial \nu}-\frac{\partial u_{i-}}{\partial \nu}\right), \quad \text{on}~ \Gamma, \\
\frac{\partial v_{i-}}{\partial \nu} - \rho \frac{\partial v_{e+}}{\partial \nu} =\frac{d}{ds}
\left( h(s) \frac{d}{ds}\left(u_{i-}-\rho u_{+}\right)\right)+ h(s)\left(k_i^2u_{i-}-k^2\rho u_{+}\right), \quad \text{on}~ \Gamma, \\
\Lim{\|{\bf x}\|\rightarrow \infty} \|{\bf x}\|^{1/2}\left(\frac{\partial v_e}{\partial r} - ikv_e\right) = 0.
\end{cases} 
\end{equation}
\end{theorem}

In the present paper, after discretizing the curve at $N=70 L k/(2\pi)$ points, where $L$ is 
the length of $\Gamma$, we apply the preceding theorem to each mode in the expansion of the 
perturbation $h(s)$, to obtain
the columns of the discretized version of $J^{(T)}_{k,\theta} \, \delta\Gamma$. From this,
we compute $(\mathcal{J}^{(T)}_k)^*$ and the steps
$\Gamma^{(j)}_{GN}$ and $\Gamma^{(j)}_{SD}$.

To summarize, the recursive linearization algorithm proceeds as outlined above, going from low
to high frequency. We refer the reader to 
\cite{borges2020inverse, borges2022multifrequency} for a more detailed discussion of 
obstacle scattering and also of the ``energy landscape" as a function of frequency. 
It is shown there, empirically, that the higher the frequency, the narrower the
basin of attraction for the nonlinear iteration and the more important it is to have 
a good initial guess.
The various parameters introduced here  ($N_{\gamma,k}, \epsilon_f$) 
are all aimed at regularizing this process in order
to obtain the best, bandlimited approximation of the unknown curve.

\section{The inverse medium problem}\label{s:vol_problem}

In this section, we briefly review the inverse medium solver of 
\cite{borges2017high}, beginning with a summary of the method used for solving the 
forward problem \eqref{eq:volume_problem} to obtain the operator
\eqref{eq:fwd_vol_op}.

Since we will be solving the same partial differential equation with many angles of
incidence, it is most efficient to use a fast direct solver. There have been 
a number of such methods proposed over the last decade (see, for example,
\cite{ambikasaran2016fast, ambikasaran2013mathcal, borm2003hierarchical, borm2003introduction, chandrasekaran2007fast, chen2002fast, corona2015n, coulier2017inverse, hackbusch2001introduction, ho2012fast, martinsson2013direct, xia2010fast, zepeda2016fast, gopal2020accelerated}). In two dimensions
(and at high frequency), most of these schemes require  $\mathcal{O}(N^{3/2})$ work to factor
the relevant linear system, where $N$ denotes the total number of degrees of freedom in the 
discretization. After factorization, the cost scales linear with $N$ for each new incident
direction (which defines the right-hand side in  \eqref{eq:volume_problem}).
As in \cite{borges2017high}, we have chosen to use
the Hierarchical Poincar\'{e}-Steklov (HPS) method of \cite{gillman2015spectrally}, since it is 
both high order accurate and efficient (that is, the constant implicit in the $O(N)$ notation
is small). 

The method begins by 
covering the domain $\Omega$ with a quad-tree data structure, designed to 
resolve $q(x)$ and to ensure a sufficient
number of points per wavelength in the discretization. In each leaf node of the tree, 
a $16\times 16$ grid of Chebyshev nodes is used to sample $q(x)$ and to 
discretize $u^{\emph{scat}}$.
Loosely speaking, the HPS method consists of three steps: 
\begin{itemize}
\item construct solution operators on leaf nodes,
\item merge solution operators on leaf nodes to construct solution operators on
their ``parent" boxes. 
\item Continue this process recursively until the solution operator for the entire
domain is available.
\end{itemize}

The actual algorithm is more complicated, and involves both an upward pass from the finest level
to the root node (a single box containing $\Omega$) and a downward pass from the root node back
to the finest level. Since the method is now well established and we use it in its standard
form, we leave a detailed description to the original paper 
\cite{gillman2015spectrally}.

\subsection{Regularization and recursive linearization}

We turn now to the inverse problem \eqref{eq:volume_inv_problem}, and first consider 
the single frequency version:
\begin{equation}
\tilde{q}=\arg \min_{q} \sum_{j=1}^{N_d} \| \ub_{k,\theta_j}^{\emph{meas}} -\mathcal{F}^{(V)}_{k,\theta_j}(q)\|. \label{eq:volume_single_inv_problem}
\end{equation}

As with the inverse obstacle problem, \eqref{eq:volume_single_inv_problem} is nonlinear, 
nonconvex, and ill-posed (without some sort of regularization). 
We will make use of the solver developed in \cite{borges2017high} which relies on
the Gauss-Newton method as a nonlinear iteration. 
Given a guess $q^{(j-1)}$, for $j=1,2,\dots$, we update the solution by letting 
$q^{(j)}=q^{(j-1)}+\delta q$ and linearize the forward 
volumetric scattering operator \eqref{eq:fwd_vol_op} to obtain the least squares problem
\begin{equation}\label{eq:volume_gn_step}
\mathcal{J}^{(V)}_k \delta q = \ub_{k}^{\emph{meas}} - \mathcal{F}_k^{(V)}(q^{(j-1)}),
\end{equation}
where $\mathcal{J}^{(V)}_k=\left(\mathcal{J}^{(V)}_{k,\theta_1};\cdots;\mathcal{J}^{(V)}_{k,\theta_{N_d}}\right)$, $\mathcal{F}^{(V)}_k=\left(\mathcal{F}^{(V)}_{k,\theta_1};\cdots;\mathcal{F}^{(V)}_{k,\theta_{N_d}}\right)$, and $\mathcal{J}^{(V)}_{k,\theta}$ is the Fr\'echet derivative of $\mathcal{F}_{k,\theta}^{(V)}$. 

The following theorem shows that the action of 
$\mathcal{J}_{k,\theta}^{(T)}$ on a {\em known} perturbation $\delta q$ can be obtained
by solving a modified scattering problem.

\begin{theorem}  (\cite{borges2017high,colton2019inverse}) 
Let $u_{k,\theta}=u_{k,\theta}^\emph{inc}+u_{k,\theta}^\emph{scat}$ denote the total field 
obtained by solving the volumetric scattering problem 
\eqref{eq:volume_problem}, where $u^\emph{inc}(\xb)= e^{ik \xb\cdot\theta}$. 
Then the forward scattering operator $\mathcal{F}^{(V)}_{k,\theta}(q)$ is Fr\'echet differentiable.
Denoting by $J_{k,\theta}^{(V)}$ the Fr\'echet derivative of 
$\mathcal{F}^{(V)}_{k,\theta}(q)$, 
let $\delta q$ be a given perturbation of the background $q(x)$. Then
the product $J_{k,\theta}^{(V)}\delta q$ is given
by the solution $v$ to the following PDE evaluated at the receiver locations:
\begin{equation}\label{eq:der_volume}
\begin{cases}
\Delta v+k^2(1+q)v = -k^2 \delta q u_{k,\theta},  \quad \text{in} ~\mathbb{R}^2 \\
\Lim{\|{\bf x}\|\rightarrow \infty} \|{\bf x}\|^{1/2}\left(\frac{\partial w}{\partial r} - ikw\right) = 0.
\end{cases} 
\end{equation}
\end{theorem}

The adjoint operator  $(J_{k,\theta}^{(V)})^*$ can be applied to a function defined at the receiver
locations in a similar fashion.

\begin{theorem}  (\cite{borges2017high}) 
Let $f(\psi)$ denote a smooth function on the circle $C$ of radius $R$ where the sensors
are located and let $\xi(f,C)$ denote the corresponding singular charge distribution on $C$ 
with charge density $f$, viewed as a generalized function in the plane.
Let $\btheta$ denote the angle of incidence of an
incoming field $u_{k,\theta}^\emph{inc}$ and let 
$q_0$ denote a known inhomogeneity in $\Omega$. Then the product
$(J_{k,\theta}^{(V)})^* \, f$ 
is given by
\[ 
(J_{k,\theta}^{(V)})^* \, f = \overline{u_{k,\theta}} w,
\]
where $u_{k,\theta}$ denotes the total field 
satisfying \eqref{eq:volume_problem} with $q=q_0$, and $w$ satisfies
\begin{equation}
\begin{cases}
\Delta w+k^2(1+q_0)w = -k^2 \xi(f,C),  \quad \text{in} ~\mathbb{R}^2 \\
\Lim{\|{\bf x}\|\rightarrow \infty} \|{\bf x}\|^{1/2}\left(\frac{\partial v}{\partial r} - ikv\right) = 0.
\end{cases} 
\end{equation}
\end{theorem}

The preceding theorems permit the solution of the 
\eqref{eq:volume_gn_step}
iteratively, using conjugate gradient iteration on the normal equations or 
LSQR iteration \cite{paige1982lsqr, barrett1994templates},
rather
than with a direct solver as in \eqref{eq:gn_step}. The number of such steps
(indexed by $j$ in equation \eqref{eq:volume_gn_step}) is controlled by a
stopping criterion. We halt when either
(a) the total number of allowed iterations $N_{it}$ has been reached, 
(b) the value of the relative residual $\|{\bf u}^\emph{meas}-\mathcal{F}_k^{(V)}(q^{(j)})\|/\|u^\emph{meas}\|<\epsilon_{res}$, 
(c) the relative size of the update $\|\delta q\|/\|q\|<\epsilon_{\delta q}$, or
(d) the residual has increased on the last step.

Recursive linearization \cite{chen1995recursive, chen1997inverse} is applied to the inverse medium problem in the same manner as
for the obstacle scattering problem, following the method
described in greater detail in \cite{borges2017high}. 
We march from low to high
frequency, with the solution obtained at frequency $k_m$ used as the initial guess 
for the solution at frequency $k_{m+1}$. By itself, this is not sufficient to overcome 
the ill-posedness of the inverse medium problem. At a given frequency, as discussed
above, the far field signature of features that are subwavelength in size is 
exponentially decaying (the Heisenberg Principle for waves). 
Thus, we bandlimit both $q(x,y)$ and the update $\delta q(x,y)$ using the 
representation
\begin{equation} \label{eq:domain_volume}
\left[ \genfrac{}{}{0pt}{}{q(x,y)}{\delta q(x,y)} \right] =
\begin{cases}
\sum\limits_{\substack{m,n=1\\ m+n\leq N_m }}^{N_m} 
\left[ \genfrac{}{}{0pt}{}{q_{m,n}}{\delta q_{m,n}} \right]
\sin(m(x+\pi/2)) \sin(n(y+\pi/2)),\quad (x,y)\in \left[-\pi/2,\pi/2\right]^2, \\
0,\quad (x,y)\notin\left[-\pi/2,\pi/2\right]^2,
\end{cases}
\end{equation}
where $N_m$ is an integer multiple of the wavenumber $k_m$. 
This imposes both compact support and regularizes the inverse problem,
controlling the condition number of the linear least squares
problem to be solved at each frequency.

\begin{remark}
At low frequencies, one can simply build the discrete version of
$J_{k,\theta}^{(V)}$ in the sine series basis and solve the least squares problem
directly (using QR factorization). At
higher frequencies, we switch to an iterative method such as 
conjugate gradient on the normal equations 
or LSQR, as noted above.
\end{remark}

\section{Numerical results}\label{s:num_res}

The stability and effectiveness of recursive linearization within the inverse obstacle solver and the inverse medium solver are sensitive to the amount of measured data available at each frequency, the contrast of the problem $k^2/k_{i}^2$, and the complexity of the domain. In this section, we explore the impact of these parameters through several numerical examples.

For each example, unless stated otherwise, we assume that scattered field measurements are made for $M$ frequencies, $k_{\ell} = 1 + (\ell-1) \delta k$, $\ell = 1,2,\ldots M$, with 
$\delta k = 0.25$, $M=117$, so that $k=30$ is the maximum frequency for which data is available. 
For each frequency, the data is obtained by solving the transmission problem where the obstacle is discretized with $100$ points per wavelength.
Scattered field measurements are made for $N_{d}$ incident waves and at $N_{r}$ receiver locations for each incident wave, with equispaced angles of incidence $\theta_{j} = (\cos{(2\pi j /N_{d}), \sin{(2\pi j/N_d))}}$, $j=1,2,\ldots N_{d}$, and equispaced receiver locations $r_{m} = 10(\cos{(2\pi m/N_{r})}, \sin{(2\pi m/N_{r})})$. For most examples, we assume $N_{d} = N_{r} = \lfloor 10k \rfloor$ at wavenumber $k$ so as to ensure that the measured data is resolved as a function of both the incidence angle and the spacing of the receiver locations. Obtaining such highly-resolved full aperture data tends to be difficult in practice. However, this idealized environment makes it feasible to study the sensitivity of the inverse solvers with respect to the contrast and complexity of the domain without a simultaneous concern about the resolution of the scattered field measurements.

For the inverse obstacle solver, the update at frequency $k$ is represented by the trigonometric series~\eqref{eq:representation_step_tranmission} with $N_{\gamma} = \lfloor 3\max{(k,k_{i})} \rfloor$ modes, while the update $\delta q$ for the inverse medium solver is parameterized by a sine series~\eqref{eq:domain_volume} with bandlimit $N_{m} = \lfloor 2k \rfloor$. The stopping criteria for both the inverse solvers were based on three considerations -- maximum number of iterations $\nitmax$, size of the relative residual $\varepsilon_{r}$, and the size of the update $\varepsilon_{u}=\|\overline{\delta \gamma}\|_2$ for the obstacle problem, and $\varepsilon_{u} = \|\overline{\delta q}\|_2/\|q\|_2$ for the medium problem. Note that we use the relative size of the update for the volume solver to account for the fact that the contrast is unknown for the inverse medium problem. For the inverse obstacle problem, we use $\varepsilon_{r} = \varepsilon_{u} = 10^{-5}$, and $\nitmax = 50$, while for the inverse medium problem, we use $\varepsilon_{r} = \varepsilon_{u} = 10^{-3}$,  and $\nitmax = 50$.
The stopping criteria have been chosen to ensure that the final reconstructions using both the solvers have converged and are not impacted by the specific choices made.

We estimate the error in the reconstruction obtained using the inverse obstacle solver by an estimate of the earth mover's distance. For this, suppose that the true obstacle boundary, and the reconstructed curve are approximated by a polygon. Let $\delta A$ denote the area of the set difference between these two polygons, and let $A$ denote the area of the true obstacle boundary. Then, $\varepsilon_{\Gamma} = \delta A/A$ is used as a measure of the relative error for the inverse obstacle problem. For the inverse medium problem, if $q$ denotes the true medium, then $\varepsilon_{q} = \|q_{r}(k) - q\|_{2}/\|q\|_{2}$ is the relative $L^2$ error in the reconstruction of the medium, where
$q_{r}(k)$ is the reconstructed medium at frequency $k$. We also compare $\varepsilon_{q}$ to the relative error of the best approximation of $q$ in the bandlimited basis given by $\varepsilon_{q_{b}} = \|q_{b}(k) - q\|_{2}/\|q\|_{2}$, where $q_{b}(k)$ is the best approximation (in an $L^2$ sense) to $q$ in the sine basis with bandlimit $N_{m}$. 

Each of the examples below is intended to 
highlight a specific feature of the inverse problem, ordered as in 
Table \ref{tab:examples}. 
\begin{table}[htbp]
	\begin{center}
		\begin{tabular}{|c | l | c|}\hline
			Example &  Description & Figures \\ \hline\hline
			1 & Impact of contrast & 1-3 \\
			2 & Effect of using limited data & 4-6\\
			3 & Reconstruction of a trapping domain & 7,8\\
			4 & Reconstruction of an obstacle with multiple components & 9\\\hline
		\end{tabular}
	\end{center}
\caption{List of numerical examples with respective results.}\label{tab:examples}
\end{table}

\subsection{Contrast}
 \label{sec:num-contrast}
\begin{figure}[h!]
\center
 \includegraphics[width=0.7\textwidth]{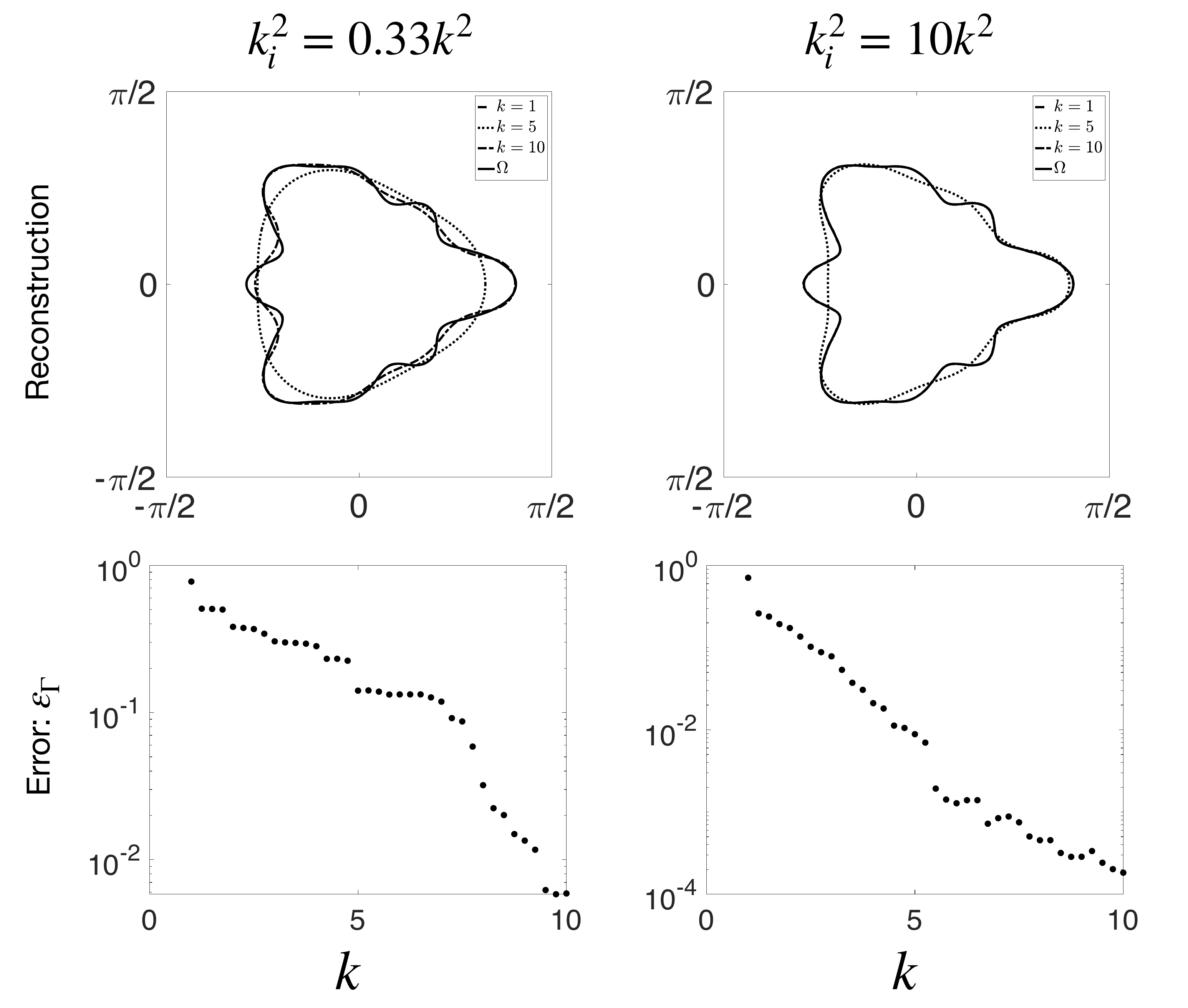}
 \caption{{\bf Effect of contrast on inverse obstacle solver:} Reconstructions of a star-shaped plane using the inverse obstacle solver and the relative error in its reconstruction $\varepsilon_{\Gamma}$.}
  \label{fig:ssp-contrast}
\end{figure}

In this section, we explore the impact of the contrast, $\eta = k_{i}^2/ k^2$ on the reconstructions obtained using the inverse obstacle and the inverse medium solvers. Consider the reconstruction of a star shaped ``glider" using the inverse obstacle solver with $\eta = 0.33$ and $\eta = 10$. The boundary of the glider is parameterized by $\gamma:[0,2\pi]\rightarrow\mathbb{R}^2$ with $\gamma(t)=r(t)(\cos(t),\sin(t))$, where $r(t)=0.9(1+\sum_{j=1}^8c_j \cos(jt)$), with $c_3=0.2$, $c_4=0.02$, $c_6=0.1$, $c_8=0.1$, and $c_j=0$, for all other values of $j$. In Figure~\ref{fig:ssp-contrast}, we plot the reconstructions at $k=1,5$, and $10$. We also plot the error $\varepsilon_{\Gamma}(k)$ corresponding to the reconstruction at wavenumber $k$.
The contrast of the object does not seem to have a significant impact on the quality of the reconstruction, and the shape of the obstacle can be recovered in a robust manner independent of the contrast of the problem.  

\begin{figure}[h!]
\center
 \includegraphics[width=0.9\textwidth]{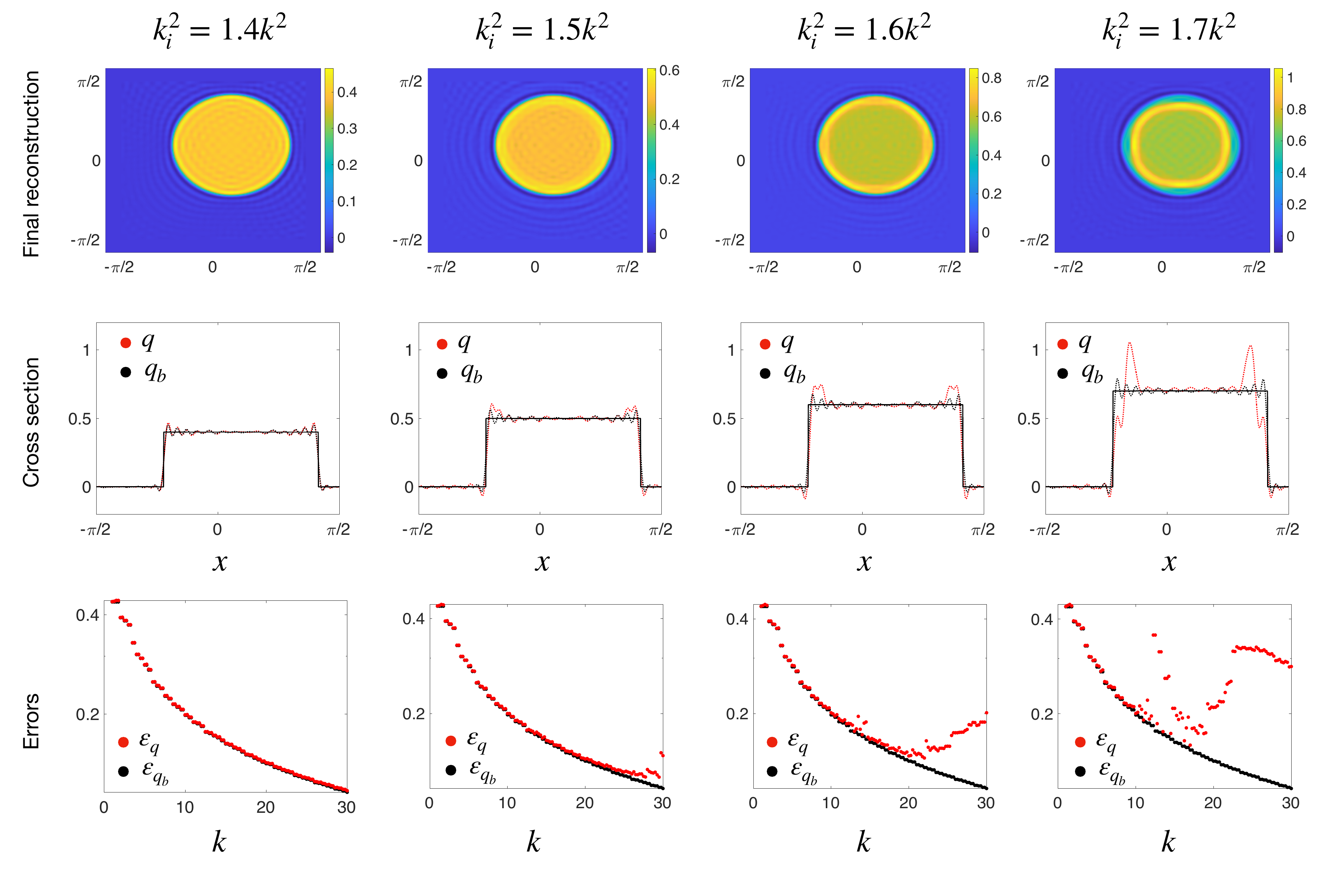}
 \caption{{\bf Effect of contrast on inverse medium solver:} Reconstructions of the circle and its relative error using the inverse medium solver. For each of the contrasts, we present the reconstruction at $k=30$, the cross section of this reconstruction through the line $(x,0.3)$, $x\in\left[-\pi/2,\pi/2\right]$ and the error of the reconstruction in the the first, second and third rows, respectively.}
  \label{fig:ex2b_rec}
\end{figure}

For the inverse medium solver, it is well known that the problem becomes easier as $\eta \to 1$, where the Born approximation leads to a more and more
accurate solution. As we increase $\eta$, the problem becomes increasingly nonlinear -- and waves interact with the inhomogeneity in more and more
complicated ways. Consider the reconstruction of a unit circle centered at $(0.3,0.3)$ with $\eta = 1.4, 1.5, 1.6,$ and $1.7$.
In Figure~\ref{fig:ex2b_rec}, we plot the final reconstructions, a cross-section of the reconstruction through the line $(x,0.3)$, and the relative errors $\varepsilon_{q}(k)$ and $\varepsilon_{q_b}(k)$ for all four cases. 
As the contrast is increased, note that the Gibbs oscillations near the boundary 
of the circle 
increase in magnitude. Moreover, the relative error $\varepsilon_{q}$ 
closely follows the 
relative error for the best reconstruction $\varepsilon_{q_{b}}$ for low contrasts, indicating that the inverse medium solver is optimal in the bandlimited basis. However, as the contrast increases, the relative error $\varepsilon_{q}$ increasingly deviates from $\varepsilon_{q_{b}}$, and the deviation begins at an earlier frequency for higher values of contrast.

\begin{figure}[h!]
 \center
\includegraphics[width=0.6\textwidth]{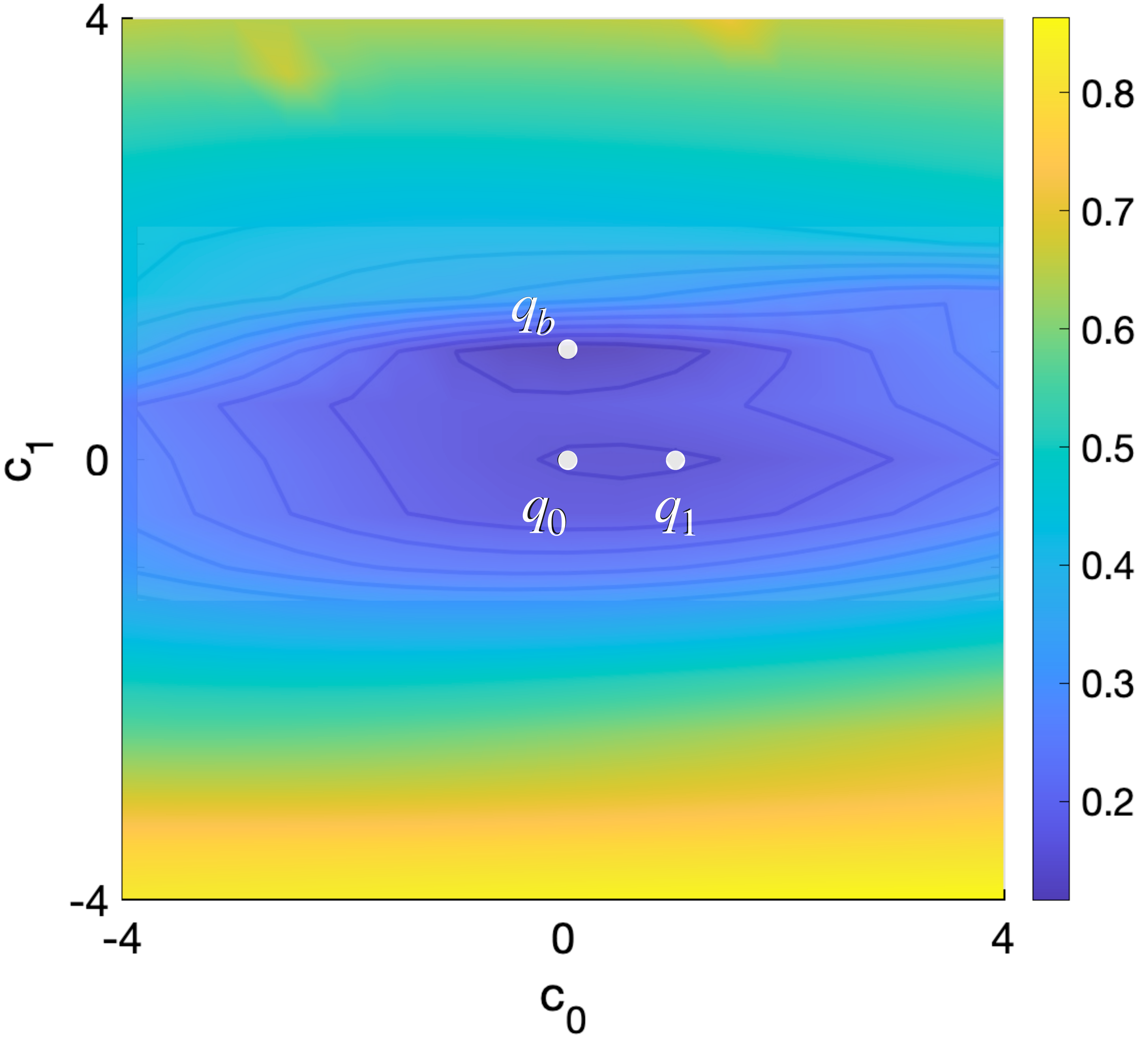}
\caption{{\bf Inverse medium problem landscape:} Contour plot of the objective function for $k=13.5$ evaluated at the plane $\mathcal{U}$ (see
\eqref{c01}). ${\bf q}_{0}$ is the initial guess, ${\bf q}_{1}$ is the solution obtained using Gauss-Newton, and ${\bf q}_{b}$ is the best solution in the bandlimited basis.}
\label{fig:local_convexity_set}
\end{figure}

This behavior can be attributed to the failure of the initial guess in 
the frequency marching process to lie within the 
basin of attraction of the {\em best} approximation.
More precisely, suppose that when $k \to k + \delta k$, 
the bandlimit of the sine series increases by $1$, i.e. $N_{m} \to N_{m}+1$. 
Following eq.~\ref{eq:domain_volume}, the new coefficients $q_{m,n} $, for $m+n = N_{m}+1$ are initialized to $0$. The local basin of attraction
of the best reconstruction, however, is approximately $O(1/k)$ in size. 
Since $N_{m} = O(k)$, for $C^{1}$ media, the Fourier coefficients decay faster than $O(1/|N_{m}|^2) = O(1/k^2)$ as $k \to \infty$. This implies that an initialization of the new coefficients to $0$ at $k+\delta k$ would roughly be $O(1/k^2)$ away from the best reconstruction in the bandlimited basis, and thus 
would likely remain in the local basin of attraction as we march in frequency. However, when recovering a discontinuous function, as is the case for our penetrable obstacle, the Fourier coefficients of the medium decay as $O(1/k)$ 
and the initialization to zero could well lie outside the 
basin of attraction. 
We can illustrate this with a simple example.
Consider the inverse problem with contrast $k_{i}^2 = 1.7 k^2$. 
Let ${\bf q}_0$ denote the initial guess at $k=13.5$, and let ${\bf q}_1$ denote the reconstruction obtained at using  Gauss-Newton. Let ${\bf q}_{b}$ denote the vector with the coefficients of the best sine series approximating the circle with the same bandlimit as ${\bf q}_0$, and ${\bf q}_{1}$. 
Let 
\begin{equation}
\label{c01}
\mathcal{U} = \left\{{\bf q}_0+c_0({\bf q}_b - {\bf q}_0)+c_1({\bf q}_1 - {\bf q}_0) \vert c_0,c_1\in\mathbb{R}\right\} .
\end{equation} 
In Figure~\ref{fig:local_convexity_set}, we plot the objective function $\|\mathcal{F}_k^{(V)}({\bf q})-u^\emph{meas}\|_{{\bf q}\in \mathcal{U}}$ in a neighborhood of ${\bf q}_0$. 
The figure shows that the initial guess ${\bf q}_{0}$ lies outside of the local set of convexity which includes the best solution ${\bf q}_{b}$, thus the reconstruction obtained using the Gauss-Newton approach ends up deviating from the best solution in the bandlimited basis. 

\subsection{Effect of limited data}

In our next example, we compare the performance of the inverse obstacle and inverse medium solvers on complicated star-shaped domains when the scattered field measurements are not completely resolved as a function of the incidence angle or receiver locations at high frequencies. The $(x,y)$ coordinates describing the boundary of the curve are given by a $35$ term Fourier series, i.e. 
\begin{equation}
\begin{bmatrix}
x(t) \\
y(t)
\end{bmatrix} = \sum_{n=-35}^{35} \begin{bmatrix}
\hat{x}_{n} \\
\hat{y}_{n} 
\end{bmatrix} e^{int} \, ,\quad t\in [0,2\pi) \, .
\end{equation}
Measurements are made for $k_{i} = 0.9 k^2$ for two sets of receiver locations and incident fields: one with $(N_{r},N_{d}) = (\lfloor 10k\rfloor, \lfloor 10k\rfloor)$, and one with $(N_{r},N_{d}) = (200,10)$. Note that the measurements are still full-aperture, however for the second case, the sensor measurements are not resolved for $k\geq 5$.  
The contrast is chosen to be low in order to isolate the effect of limited data on the reconstructions.

\begin{figure}[h!]
 \center
\includegraphics[width=\textwidth]{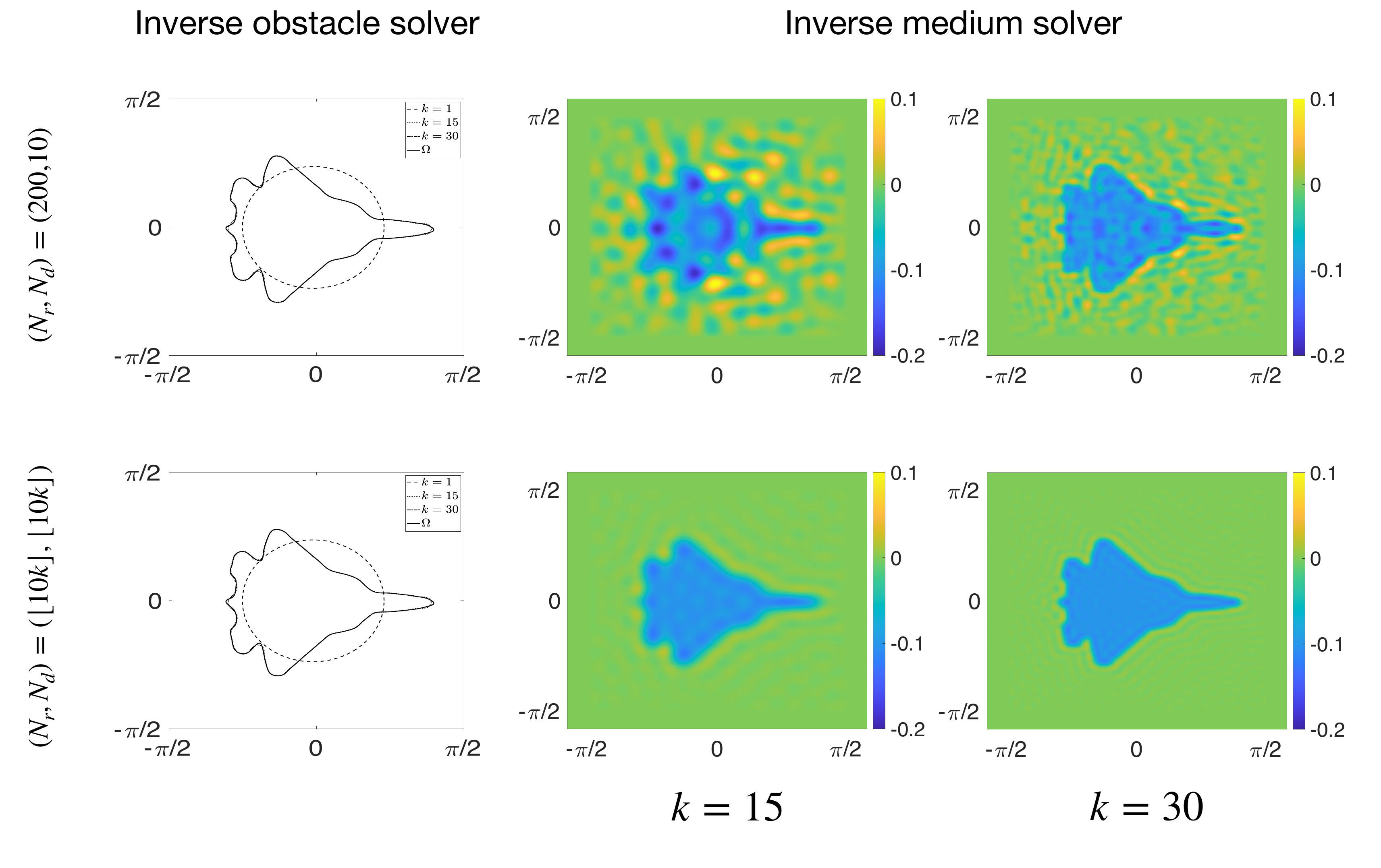}
\caption{{\bf Reconstructions for the limited data experiment}: Reconstruction of a star-like plane using the inverse obstacle solver at frequencies $k=1,15$, and $30$ (left), the inverse medium solver at $k=15$ (middle), or the inverse medium solver at $k=30$ (right). The top row results correspond to 
underresolved scattered field measurements at high frequencies, while the bottom row results correspond to fully resolved scattered field measurements. }
\label{fig:ex3_plane}
\end{figure}

\begin{figure}[h!]
 \center
\includegraphics[width=\textwidth]{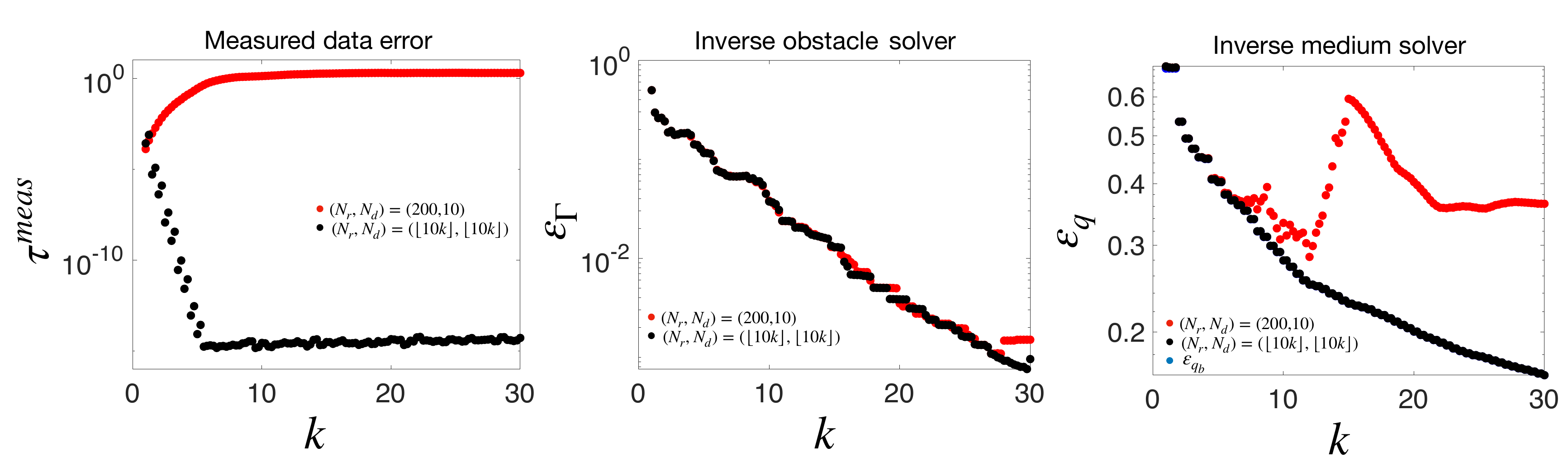}
\caption{ {\bf Errors in reconstruction and measurement data for the limited data experiment:} (left) Resolution of measurement data $\tau^{\textrm{meas}}$, (middle) error in reconstruction for the inverse obstacle solver $\varepsilon_{\Gamma}$, and (right) error in reconstruction for the inverse medium solver $\varepsilon_{q}$ along with the error in the best reconstruction in the bandlimited basis $\varepsilon_{q_{b}}$, all as a function of frequency.}
\label{fig:err_limited_data}
\end{figure}

In Figure~\ref{fig:ex3_plane}, we plot the reconstructions obtained using the inverse obstacle solver at $k=1$, $15$, and $30$, and the inverse medium solver at $k=15$ and $30$, and in Figure~\ref{fig:err_limited_data}, we plot the error in reconstruction $\varepsilon_{\Gamma}(k)$ for the inverse obstacle solver, and the error in reconstruction for the inverse medium solver $\varepsilon_{q}(k)$, along with the error for the best solution in the bandlimited basis $\varepsilon_{q_{b}}$.
To illustrate the lack of resolution of the measured data, 
we also plot 
{\small
\begin{equation}
\tau^{\textrm{meas}} = \frac{1}{\max_{m,n} |\hat{u}_{m,n}|} \left( \sum_{m=-\frac{N_{r}}{2}+1}^{\frac{N_{r}}{2}} \left( |\hat{u}_{m,-N_{d}/2+1}| + |\hat{u}_{m,N_{d}/2}| \right) + \sum_{n=-\frac{N_{d}}{2}+1}^{\frac{N_{d}}{2}} \left( |\hat{u}_{-N_{r}/2+1,n}| + |\hat{u}_{N_{d}/2,n}| \right) \right) \, ,
\end{equation}}
where $\hat{u}_{m,n}$ are the Fourier series coefficients of the measured data, i.e.
\begin{equation}
u^{\textrm{meas}}(\theta,\phi) = \sum_{m=-\frac{N_{r}}{2}+1}^{\frac{N_{r}}{2}} \sum_{n=-\frac{N_{d}}{2}+1}^{\frac{N_{d}}{2}} \hat{u}_{m,n} e^{i(m\theta + n \phi)} \, .
\end{equation}

Both the inverse medium solver, and inverse obstacle solver recover the obstacle to high fidelity when sufficient data is available to resolve the scattered field measurements. However, the inverse medium solver deviates from the best solution in the bandlimited basis when the data is insufficient to resolve the scattered field (approximately when $\tau^{\textrm{meas}} > 0.1$), while the inverse obstacle solver robustly recovers the shape even with unresolved scattered field measurements. 
The inverse obstacle solver tends to perform better than the inverse medium solver with fewer measurements of the scattered field, since at wavenumber $k$, the obstacle is described by $O(k)$ parameters while the medium is described by $O(k^2)$ parameters.

\begin{figure}[h!]
 \center
\includegraphics[width=0.7\textwidth]{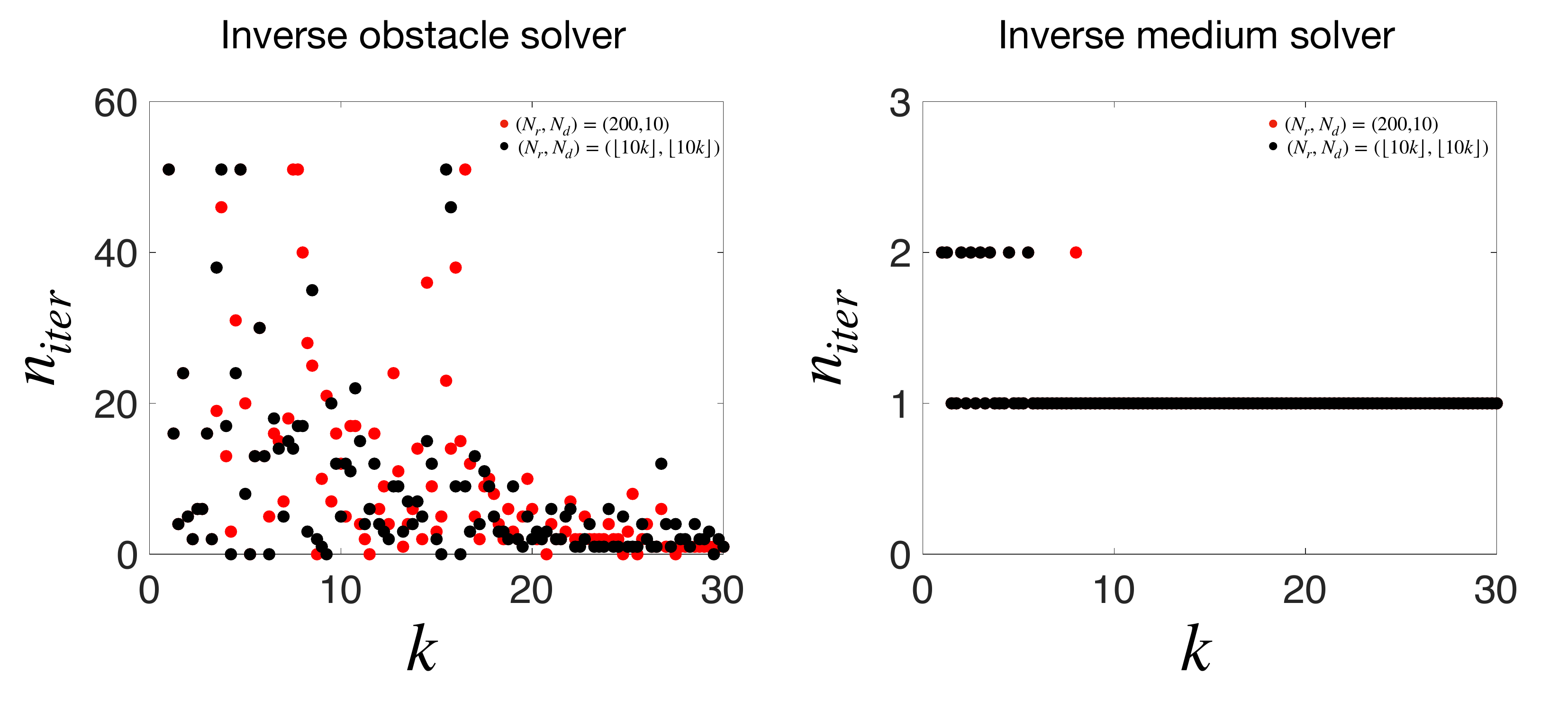}
\caption{{\bf Iteration count for inverse solvers in limited data experiment:} Number of iterations $\niter$ in the optimization loop at frequency $k$ for the inverse obstacle solver (left), and the inverse medium solver (right).}
\label{fig:limited_iter}
\end{figure}

Finally, in Figure~\ref{fig:limited_iter}, we plot the number of iterations $\niter$ in the optimization loop at each frequency.  It is interesting to note that the inverse medium solver tends to take fewer steps when compared to the inverse obstacle solver. This behavior is independent of the shape of the obstacle being recovered, and can be attributed, perhaps, to the more systematic 
increase in resolution that is achieved via the tensor-product sine series in the volumetric case. 
That said, the computational complexity of both of these approaches is comparable, since the obstacle problem benefits from dimensionality reduction:
one needs to solve only a boundary integral equation instead of a volumetric
Lippmann-Schwinger equation.

\subsection{Cavity-like domains}
\begin{figure}[h!]
 \center
\includegraphics[width=\textwidth]{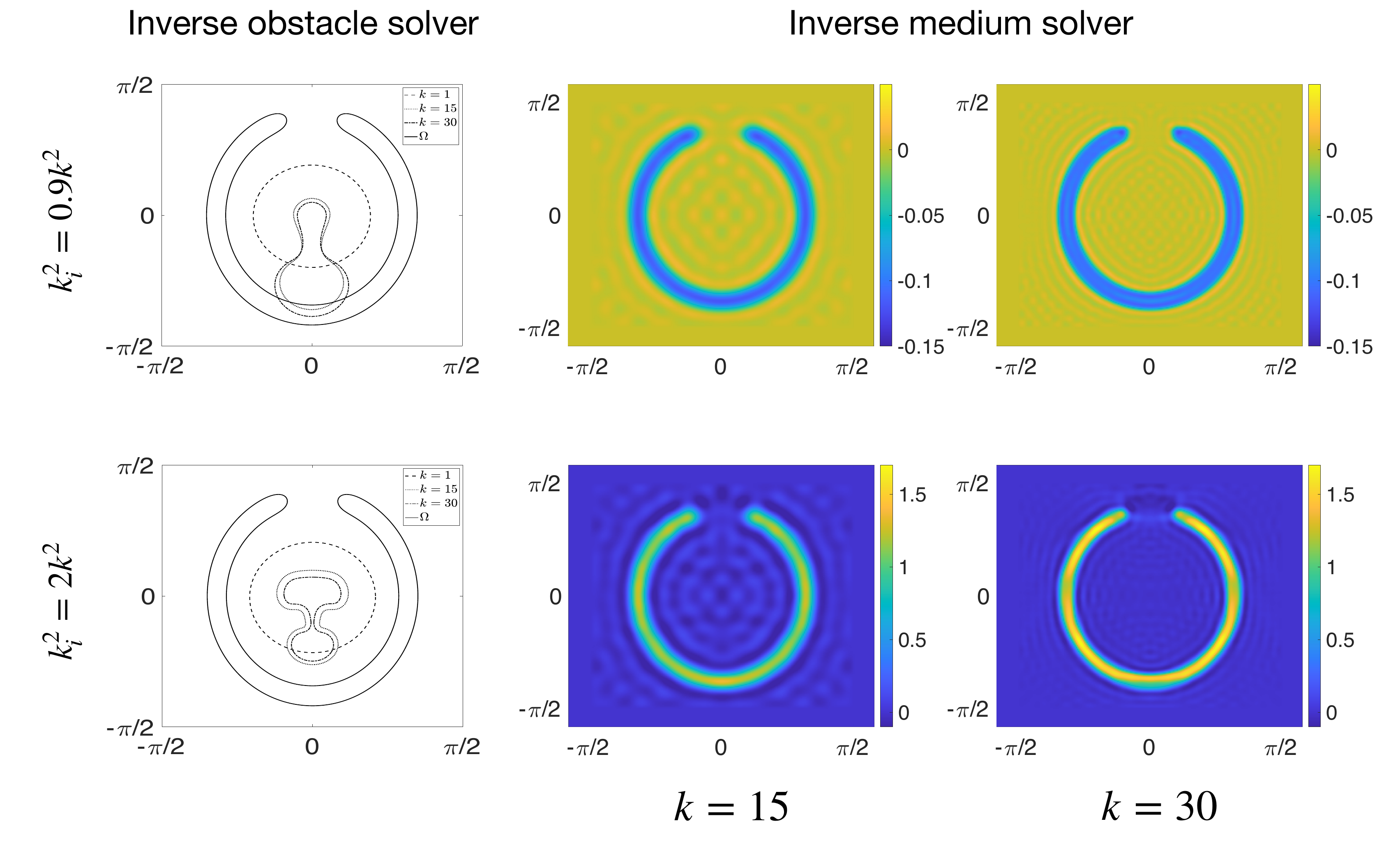}
\caption{{\bf Reconstructions of a trapping domain:} Reconstruction of a trapping cavity-like domain using the inverse obstacle solver at frequencies $k=1$, $15$, and $30$, using the inverse medium solver at $k=15$, and using the inverse medium solver at $k=30$, in the left, middle and right columns. The top row results are for $k_{i}^2 = 0.9 k^2$, and the bottom row results correspond to $k_{i}^2 = 2k^2$.}
\label{fig:cavity_recons}
\end{figure}

\begin{figure}[h!]
 \center
\includegraphics[width=0.7\textwidth]{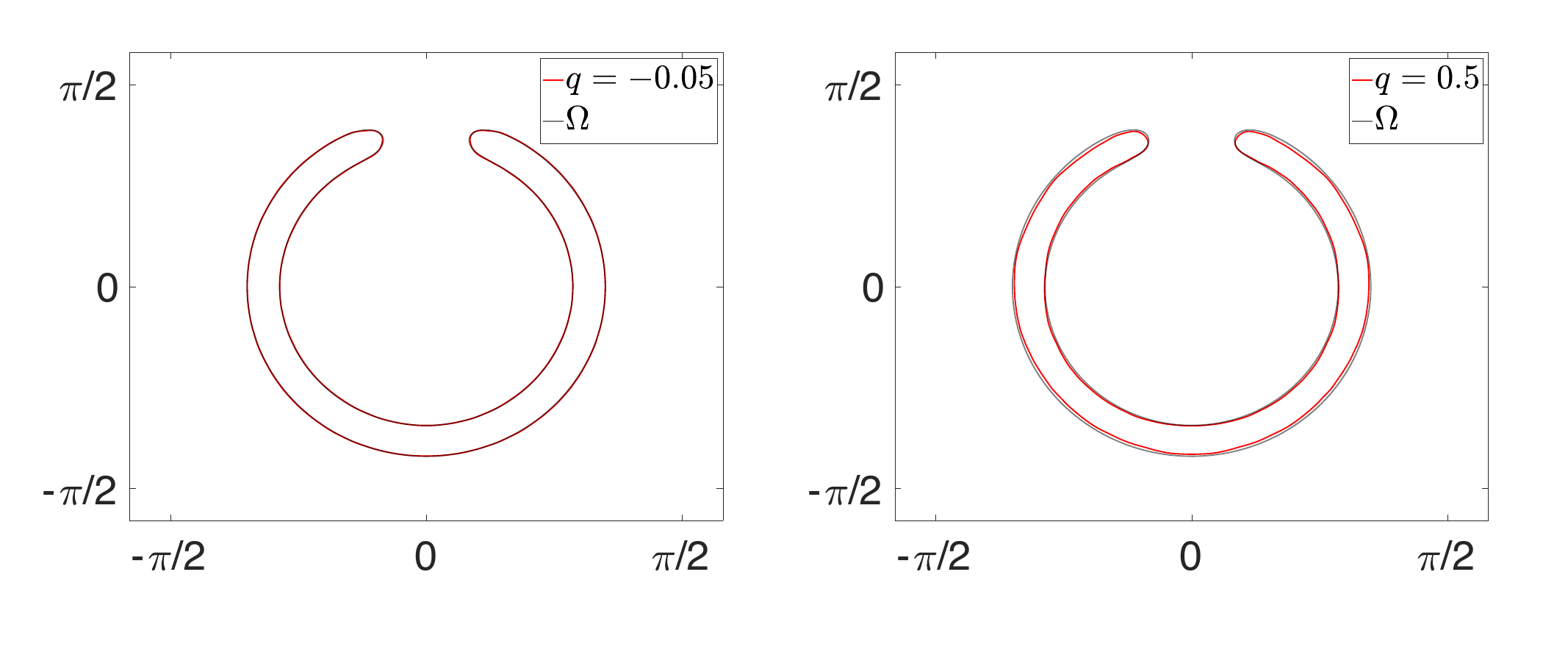}
\caption{{\bf Level sets from the inverse medium solver for the trapping domain:} (left) level set for $k_{i}^2 = 0.9 k^2$, and (right) level set for $k_{i}^2 = 2k^2$. The boundary of the curve used to generate the data is plotted in gray in both figures.}
\label{fig:cavity_levelset}
\end{figure}

In our next example, we investigate the behavior of the two solvers for (highly
nonconvex) cavities, which are challenging;
the solutions on such domains tend to have a complicated behavior as a function of frequency owing to the wave-trapping nature of the domain. The cavity used for this example is illustrated by the curve in the left-hand panels
of Figure~\ref{fig:cavity_recons}.
In that figure, we plot the reconstructions obtained using the inverse obstacle solver and the inverse medium solver at $k=15$, and $k=30$ for two different contrast values $k_{i}^2 = 0.9 k^2$ and $k_{i}^2 = 2k^2$. The results show that the inverse obstacle solver is unable to resolve the shape of the domain, while the inverse medium solver performs significantly better.
This is further illustrated in Figure~\ref{fig:cavity_levelset}, where we plot the level set corresponding to $q/2$, namely
$q=-0.05$ in the low contrast case and $q=0.5$ in the higher contrast case.
This is where one would expect the boundary to lie in a truncated Fourier series
subject to the Gibbs phenomenon.

\subsection{Multiple obstacles}
\begin{figure}[h!]
 \center
\includegraphics[width=\textwidth]{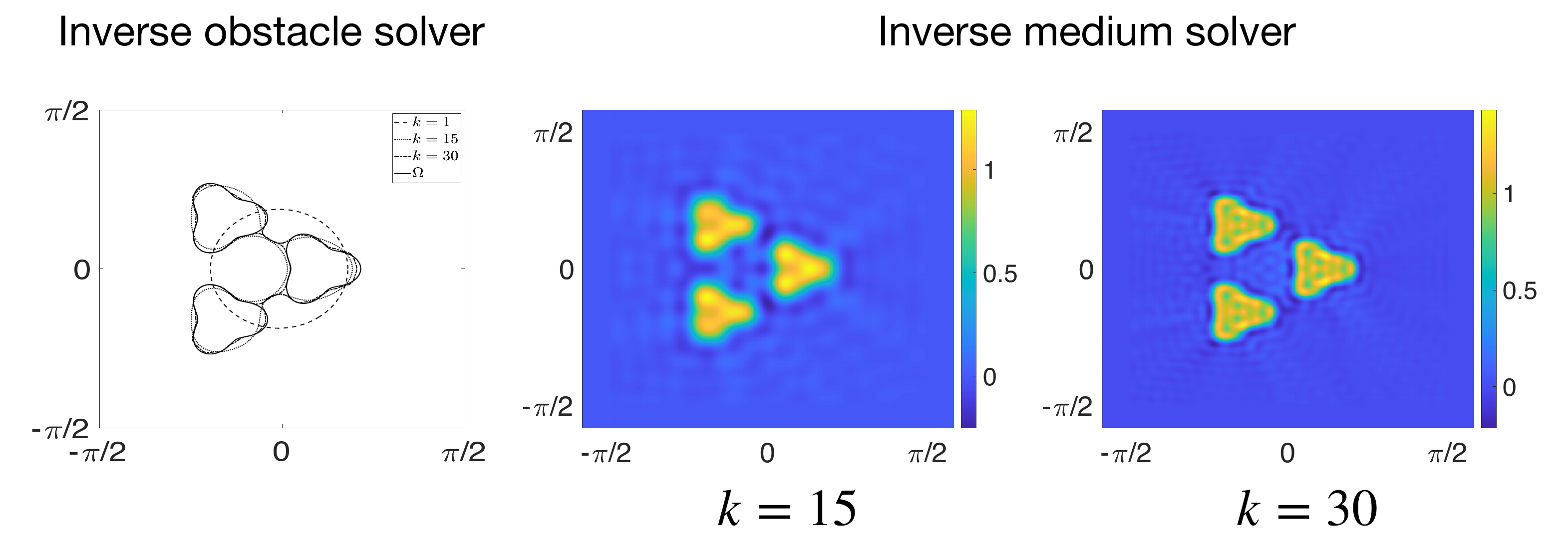}
\caption{{\bf Reconstructions of multiple scatterers: } Reconstructions of three star shaped scatterers using the inverse obstacle solver at frequencies $k=1,15$, and $30$, using the inverse medium solver at $k=15$, and using the inverse medium solver at $k=30$, in the left, middle and right columns.}
\label{fig:multiple-scat}
\end{figure}

For our final example, we consider the reconstruction of multiple, disjoint 
scatterers 
using both the inverse obstacle and inverse medium solvers. The domain consists of three identical star-shaped scatterers with $k_{i}^2 = 2k^2$ for each 
(see the left panel of Figure~\ref{fig:multiple-scat}). 
At wavenumber $k$, measurements are made for $\lfloor 10k \rfloor$ incident directions at $\lfloor  10k \rfloor$ sensor locations. While the inverse medium solver should not be impacted by the presence of multiple scatterers, 
this is a challenging problem for the inverse obstacle solver since our model
is parametrized as a single closed curve.

In Figure~\ref{fig:multiple-scat}, we plot the reconstructions obtained using both solvers at $k=15$, and $k=30$. 
Since the obstacles are separated by more than a wavelength corresponding to the smallest wavelength for which measurements are made, the inverse medium solver is robustly able to recover the boundary of the multiple obstacles.  
Perhaps surprisingly, the inverse obstacle solver also performs well, 
and is able to capture the bulk of the boundary accurately, with the three scatterers connected via thin bridges. Detecting such a feature could serve
as a signal that multiple obstacles are present, which could be used, in turn,
to modify the number of components used in the model.

 \section{Conclusions}\label{s:conclusions}
In this paper, we compare the performance of an inverse obstacle solver for transmission boundary value problems and an inverse medium solver for the reconstruction of penetrable media. For both solvers, full aperture measurements are made at multiple frequencies, and a recursive linearization based approach is used as a continuation method for solving the a sequence of single frequency 
optimization problems. 

We present several numerical examples which highlight the benefits and disadvantages of the two approaches. The inverse obstacle solver requires the determination of a simple curve, while the inverse medium solver discretizes an entire 
two dimensional volume. 
We have found that the inverse obstacle solver tends to perform better 
when the number of incidence directions and receiver locations are insufficient to resolve the scattered field measurements. This is to be expected, since
many fewer degrees of freedom are required for describing the unknown obstacle. 
On the other hand, 
the constrained optimization problem for the inverse medium solver appears
to be better posed, since the space of compactly supported perturbations 
comes equipped with a natural basis (such as a sine series), and it is 
straightforward to systematically increase the bandlimit. 
For the inverse obstacle problem, the geometry of the set of non-intersecting curves is complicated to parametrize, and constructing a constraint set which appropriately bandlimits the curve poses a significant challenge and local linearization appears to be less robust. The difference in the optimization landscape between these two problems manifests itself in two ways. First, for all obstacles, the inverse medium solver tends to require fewer optimization steps at each frequency. Second, the inverse medium solver is more robust for reconstructing complicated shapes such as cavity-like structures, where the measurement data can vary sharply across frequencies, due to the trapping nature of the domain. 

However,
the inverse medium solver is more sensitive to contrast, deviating from the best solution in the bandlimited basis as the contrast is increased. This behavior could be due to poor initialization of the Newton iteration as we 
march in frequency - a question which is currently being explored.

Given the better stability of the inverse obstacle solver when scattered field measurements are limited, and the better behavior of the inverse medium solver in complicated cavity-like domains, one can imagine using 
the inverse medium solver at low frequencies to construct a good initial guess
and to continue from that point using the inverse obstacle solver at higher frequencies. Such hybrid schemes are under investigation.

Finally, we showed that both solvers were able to handle the case of 
multiple scatterers. For the inverse medium solver, this is to be expected,
since the collection of scatterers is treated together as an unknown
function to be recovered. We were surprised that, without prior knowledge of 
the number of scatterers, the inverse obstacle solver was able to recover 
the geometry with high fidelity, stitching together 
the scatterers with very thin bridges. We suspect that the presence of 
such thin bridges could be used as a monitor to detect the presence of 
multiple scatterers and the reconstruction could be further refined 
by adjusting the number of boundary curves during the optimization process.
This is also an area of ongoing research.

\section*{Acknowledgments}
The work of C.~Borges was supported in part by the Office of Naval Research under award number N00014-21-1-2389. The authors would like to thank Jeremy Hoskins and Travis Askham for many useful discussions. 

\bibliography{./biblio_comparison}

\end{document}